\documentclass[dvips,aos,preprint]{imsart}

\RequirePackage[OT1]{fontenc}

\usepackage{amsmath,amssymb,amsthm, latexsym,amsfonts}
\usepackage[dvips]{graphics,epsfig}

\newcommand{\N}{\mathbb{N}}

\newcommand{\BB}{\mathcal{B}}

\numberwithin{equation}{section}

\def\Min(#1,#2){#1\wedge #2}
\def\Max(#1,#2){#1\vee #2}

\def\rueck{\noindent\hangafter=1 \hangindent=1.3em}

\begin{document}

\begin{frontmatter}

\title{Correction note to ``Limit Theorems for Empirical Processes of Cluster Functionals''}
\runtitle{Empirical Cluster Processes}
\thankstext{T1}{We would like to thank Johan Segers
for pointing out the error discussed in this note.}

\begin{aug}
  \author{\fnms{Holger}  \snm{Drees}\corref{}\ead[label=e1]{holger.drees@uni-hamburg.de}},
  \and \author{\fnms{Holger} \snm{Rootz\'{e}n}\ead[label=e2]{rootzen@math.chalmers.se}}

  \affiliation{University of Hamburg and Chalmers and Gothenburg University}

  \address{University of Hamburg\\ Department of Mathematics,
SPST\\ Bundesstr.\ 55\\ 20146 Hamburg\\ Germany\\
          \printead{e1}}

  \address{Chalmers University \\ Department of Mathematical Statistics\\
412 96 G\"{o}teborg\\ Sweden\\ and\\
 Gothenburg University\\ Department
of
Mathematical Statistics\\  412 96 G\"{o}teborg\\ Sweden\\
          \printead{e2}}

\end{aug}



\end{frontmatter}

In Lemma 5.2 (vii) it is stated that under the conditions (B1) and (B3) the length $L(Y_n)$ of the core of a cluster  satisfies $\lim_{k\to\infty} \limsup_{n\to\infty} P\{L(Y_n)>k\}/(r_nv_n)=0$. However, in general, the first inequality in the proof of this part of the lemma is not correct, and it seems likely that the assertion does not hold under the stated conditions. Note that this part of the lemma is used only in Remark 3.7 (i); so none of the other results are affected.

The easiest way to correct the error is to replace condition (B3) with the corresponding condition for $\varphi$-mixing coefficients
$$ \varphi_{n,k} := \sup_{1\le l\le n-k-1} \sup_{B\in\BB_{n,l+k+1}^n, C\in\BB_{n,1}^l}\big| P(B)-P(B|C)\big|
$$
(with the convention $P(B|C)=P(B)$ if $P(C)=0$), i.e.\ to assume $\lim_{m\to\infty}$ $ \limsup_{n\to\infty} \varphi_{n,m}=0$. The arguments of the proof of Lemma 5.2 (vii) are rectified if $\beta_{n,k}$ is replaced with $\varphi_{n,k}$ everywhere.

However, often the following simpler condition is easier to verify:\\[1ex]
\parbox[t]{1cm}{$(\widetilde{B3})$} \parbox[t]{11.5cm}{For all $n\in\N$ and all $1\le i\le r_n$ there exists $s_n(i)\ge P(X_{n,i+1}\ne 0\mid X_{n,1}\ne 0)$ such that $s_\infty(i):=\lim_{n\to\infty} s_n(i)$ exists and $\lim_{n\to\infty} \sum_{i=1}^{r_n} s_n(i) =\sum_{i=1}^\infty s_\infty(i)<\infty$.}
\smallskip

Since, by stationarity,
\begin{eqnarray*}
\frac 1{r_nv_n} P\{L(Y_n)>k\} & \le & \frac 1{r_nv_n} \sum_{i=1}^{r_n-k} \sum_{j=i+k}^{r_n} P(X_{n,j}\ne 0|X_{n,i}\ne 0) P\{X_{n,i}\ne 0\} \\
& \le & \sum_{j=k}^{r_n} s_n(j),
\end{eqnarray*}
the assertion of Lemma 5.2 (vii) follows readily.

To check condition $(\widetilde{B3})$, typically one bounds $P(X_{n,i+1}\ne 0\mid X_{n,1}\ne 0)$ by an expression of the form $s_n(i)=b_n+c_i$ with $b_n=o(1/r_n)$ and $\sum_{i=1}^\infty c_i<\infty$. The interchangeability of the limit and the sum  is then automatically fulfilled. For example, $(\widetilde{B3})$ has been verified in Example 8.3 of Drees et al.\ (2015) for solutions to stochastic recurrence equations.

Condition $(\widetilde{B3})$ has the additional advantage that in Remark 3.7 (i) it renders condition (3.9) superfluous, i.e.\ condition (C3) is met if $(\widetilde{B3})$ and (3.8) hold. To see this, check that, for bounded functions $\phi,\psi$, using stationarity $E(g_\phi(Y_n)g_\psi(Y_n))/(r_nv_n)=Cov(g_\phi(Y_n),g_\psi(Y_n))/(r_nv_n)+O(r_nv_n)$ can be represented as
\begin{eqnarray*}
  \lefteqn{\frac 1{v_n} E\big(\phi(X_{n,1})\psi(X_{n,1})\big)}\\ & & + \sum_{k=1}^{r_n-1} \frac 1{v_n}\Big(1-\frac{k}{r_n}\Big) \big(E(\phi(X_{n,1})\psi(X_{n,k+1})+E(\psi(X_{n,1})\phi(X_{n,k+1})\big)
\end{eqnarray*}
which tends to $c(g_\phi,g_\psi)$ defined in (3.10) by our assumptions and Pratt's lemma (Pratt, 1960), because the $k$-th summand can be bounded in absolute value by
$ 2\|\phi\|_\infty\|\psi\|_\infty s_n(k)$. Moreover, using the above representation with $\phi=\psi=1_{E\setminus\{0\}}$, one immediately sees that $(\widetilde{B3})$ also implies condition (3.5).

\medskip

\noindent {\large\bf References}
   \smallskip

\parskip1.2ex plus0.2ex minus0.2ex

\rueck
Drees, H., and Rootz\'{e}n, H. (2010). Limit Theorems for Empirical Processes of Cluster Functionals. {\em Ann.\ Statist.} {\bf 38}, 2145--2186.

\rueck
Drees, H., Segers, J., and Warcho\l, M. (2015). Statistics for Tail Processes of Markov Chains. {\em Extremes} {\bf 18}, 369--402.
\smallskip

\rueck
Pratt, J.\ W.\ (1960). On interschanging limits and integrals. {\em Ann.\ Math.\ Statist.} {\bf 31}, 74--77.

\end{document}